\documentclass{amsart}

\newtheorem{theorem}{Theorem}[section]
\newtheorem{lemma}[theorem]{Lemma}
\newtheorem{proposition}[theorem]{Proposition}
\newtheorem{corollary}[theorem]{Corollary}

\theoremstyle{definition}

\newtheorem{remark}[theorem]{Remark}
\newtheorem{question}[theorem]{Question}
\newtheorem*{acknowledgment}{Acknowledgment}

\numberwithin{equation}{section}

\usepackage{amssymb, amscd, amsfonts}

\newcommand{\QQ}{\mathbb{Q}}

\newcommand  {\shI}     {\mathcal{I}}

\newcommand  {\foX}     {\mathfrak{X}}
\newcommand  {\foY}     {\mathfrak{Y}}

\newcommand  {\aff}     {{\text{aff}}}

\newcommand  {\Bs}      {\operatorname{Bs}}

\newcommand  {\dom}     {\operatorname{dom}}

\newcommand  {\lra}     {\longrightarrow}

\newcommand  {\mult}    {\operatorname{mult}}

\renewcommand{\O}       {\mathcal{O}}

\renewcommand{\P}       {\overline{P}}

\newcommand  {\Pic}     {\operatorname{Pic}}

\newcommand  {\Proj}    {\operatorname{Proj}}

\newcommand  {\ra}      {\rightarrow}

\newcommand  {\SBs}     {\operatorname{SBs}}

\newcommand  {\Spec}    {\operatorname{Spec}}

\def\mydate{\number\day\space\ifcase\month \or January\or February\or March\or 
April\or May\or\June\or\July\or
August\or September\or October\or November\or December\fi \space\number\year}

\begin{document}

\title[Canonical models of normal surfaces ]{What is missing in  canonical 
models 
for proper normal algebraic surfaces?}

\author[Stefan Schroer]{Stefan Schr\"oer}
\address{Mathematische Fakult\"at, Ruhr-Universit\"at, 
               44780 Bochum, Germany}
\curraddr{Massachusetts Institute of Technology,
Department of Mathematics, Building 2 Room 155, 
77 Massachusetts Avenue, Cambridge MA 02139-4307, USA}

\email{s.schroeer@ruhr-uni-bochum.de}

\subjclass{14C20, 14J17, 14J29}

\dedicatory{Revised version, 9 July 2001}

\begin{abstract}
Smooth surfaces have finitely generated canonical rings and   projective 
canonical models.
For normal surfaces, however, the graded ring of multicanonical sections  is 
possibly nonnoetherian,
such that the corresponding homogeneous spectrum is noncompact.
I construct a canonical
compactification   by adding finitely many non-$\QQ$-Gorenstein points at 
infinity, provided
that  each Weil divisor is numerically equivalent to a $ \QQ$-Cartier divisor.
Similar results hold for arbitrary Weil divisors instead of the canonical
class.
\end{abstract}

\maketitle

%===========================================================
\section*{Introduction}

Each proper normal algebraic surface 
$X$ comes along with the  graded ring 
$R(K_X)= \bigoplus_{n\geq 0} H^0(nK_X)$ and the    homogeneous spectrum
$P(K_X)=\Proj R(K_X)$.
If 
$X$ has   canonical singularities (that is, rational Gorenstein 
singularities), 
these are  
called the \emph{canonical ring} and the \emph{canonical model}.
It is then  a classical 
theorem, and   Mori theory reassures us, 
that the ring
$R(K_X)$ is finitely generated, such that the scheme
$P(K_X)$ is projective.

A natural question to ask: What happens for proper  surfaces
with arbitrary normal singularities? The usual approach is 
to pass to  smooth models,
but I want so see what happens on the singular surface itself.
To my knowledge, neither canonical rings nor canonical models were  
 studied from this viewpoint,
possibly because
Zariski \cite{Zariski 1962} observed that the   ring 
$R(K_X)$ might be nonnoetherian. However,  $P(K_X)$ is 
always a scheme of 
finite type \cite{Schroeer 2000},
so it makes perfect sense to analyze it from a geometric point of view.  

According to the Nagata Compactification Theorem, we can embed 
$P(K_X)$ into a proper normal scheme by adding points and curves at infinity.
 I prefer  to add just points,
because such compactifications are minimal and therefore unique. 
In this paper, we 
shall construct such
a compactification 
$P(K_X)\subset\P(K_X)$ with discrete boundary by adding finitely many 
non-$\QQ$-Gorenstein points at infinity.
I have, however, to assume that  each Weil divisor on 
$X$ is numerically equivalent to a 
$\QQ$-Cartier divisor. This  holds, for example,  for surfaces with geometric 
genus 
$p_g=0$, or for surfaces with rational singularities.
Much of this generalizes to  $P(D)$, where $D\in Z^1(X)$ is an arbitrary Weil 
divisor.

The paper is organized as follows. In Section 1, we collect some facts on the
intersection of the base loci $\bigcap \Bs(nD)$.  In Section 2, 
we determine the 
structure
 of the rational map $r_D:X\dashrightarrow P(D)$. In the next section, we use 
this rational map
 and construct a compactification $P(D)\subset \P(D)$ as an  algebraic space.
Section 4 contains some results on the multiplicities occurring in the base
loci $\Bs(nD)$. In the last section, we apply our results to the 
canonical model
$P(K_X)$.

\begin{acknowledgment}
I wish to thank Hubert Flenner for drawing my attention to the homogeneous 
spectrum
$P(K_X)$.
Furthermore, I thank the
Mathematical Department of the
Massachusetts Institute of Technology for its hospitality,
and the Deutsche Forschungsgemeinschaft for financial support.
Finally, I thank the referee for his comments, which helped to
clarify the paper.
\end{acknowledgment}

%===========================================================
\section{Stable base locus}

Throughout the paper, 
$X$ is a proper normal algebraic surface. 
In other words, a normal 2-dimensional 
scheme 
proper over an arbitrary  base  field $k$. 
Each Weil divisors
$D\in Z^1(X)$ yields a graded algebra 
$R(D)=\bigoplus_{n\geq 0}H^0(nD)$, which in turn defines  the   
homogeneous spectrum 
$P(D)=\Proj  R(D)$. Here we use  
$H^0(nD)$ as a shorthand for $H^0(X,\O_X(nD))$. 
Let me quote the following fact (\cite{Schroeer 2000},
Section 6).

\begin{proposition}
\label{algebraic}
The homogeneous spectrum 
$P(D)$ is a normal algebraic 
scheme of dimension 
$\leq 2$. Furthermore, if 
$\dim P(D)\leq 1$, then 
$P(D)$ is projective.
\end{proposition}

The interesting point  is that $P(D)$ is of finite type over $k$ (which
does not necessarily hold in higher dimensions).
By definition, there is  an affine open covering
$$
P(D)= \bigcup \Spec(R(D)_{(s)})
$$
where the union runs over all homogeneous 
$s\in R(D)$ of positive degree. Note that, if 
$X_s\subset X$ is the open subset where 
$s:\O_X\ra\O_X(nD)$ is bijective, we have 
$\Gamma(X_s,\O_X)=R(D)_{(s)}$. Thus we can write 
$$
P(D)= \bigcup X_s^\aff,
$$
where 
$X_s^\aff=\Spec\Gamma(X_s,\O_X)$ is the \emph{affine envelope}. Gluing the 
canonical morphisms 
$U\ra U^\aff$ gives a rational map 
$r_D:X\dashrightarrow P(D)$.

The \emph{base locus} 
$\Bs(nD)\subset X$ is the intersection of all curves 
$C\sim D$ (linear equivalence). Following Fujita \cite{Fujita 1983},
Definition 1.17, 
we define the \emph{stable base locus}   
$$
\SBs(D)=\bigcap_{n>0}\Bs(nD)\subset X.
$$
Obviously, the rational map 
$r_D:X\dashrightarrow P(D)$ is defined on $X-\SBs(D)$.
Let us collect some facts on stable base loci.
Clearly, 
$x\in\SBs(D)$ if the divisor germ 
$D_x\in Z^1(\O_{X,x})$ is not 
$\QQ$-Cartier. There is a partial converse as follows. 
Define 
$$
\SBs^0(D)\subset\SBs(D)
$$
to be the 0-dimensional part of the stable base locus.

\begin{proposition}
\label{not Q-Cartier}
For each  
$x\in\SBs^0(D)$, the divisor  germ 
$D_x\in Z^1(\O_{X,x})$ is not 
$\QQ$-Cartier.
\end{proposition}

\proof
Clearly 
$\SBs(D)=\SBs(nD)$  for all integers 
$n>0$, so  we may assume that 
$D$ is a curve. Seeking a contradiction, we assume that 
$D_x$ is 
$\QQ$-Cartier. Passing to a suitable multiple, we may assume that 
it is Cartier. Set $\shI=\O_X(-D)$. The blowing up 
$Y=\Proj(\bigoplus\shI^n)$   and the invertible sheaf 
$\O_Y(1)$ depend only on the linear equivalence class of $D$.
Let $\tilde{Y}$ be the normalization of the blowing-up. The induced morphism 
$f:\tilde{Y}\ra Y$ is bijective near $x$, and $y=f^{-1}(x)$ 
is an isolated base 
point for the Cartier divisor $f^{-1}(D)$. The latter contradicts the 
Fujita--Zariski Theorem (\cite{Fujita 1983}, Theorem 1.19).
\qed

\medskip
Next, we consider the 1-dimensional part of the stable base locus.

\begin{proposition}
\label{bijective support}
Suppose 
$D\in Z^1(X)$ is  a Weil divisor with 
$H^0(nD)\neq 0$ for some 
$n>0$. Let 
$E\subset X$ be a curve supported by 
$\SBs(D)$. Then the canonical map 
$H^0(D)\ra H^0(D+E)$ is bijective.
\end{proposition}

\proof
Inducting on the number of irreducible components of $E$, we may assume that 
$E$ is irreducible. Assuming that
$E$ is also reduced, we have to check that the canonical inclusion
$H^0(D)\subset H^0(D+mE)$ is bijective for all 
$m>0$. Suppose to the contrary that  
$D+mE\sim A+m'E$ for some 
$0\leq m'<m$ and a curve 
$A\subset X$ not containing $E$. Subtracting 
$m'E$, we may assume 
$D +mE\sim A$. Decompose 
$nD= B+lE$ for some integer 
$l\geq 1$ and a curve 
$B\subset X$ not containing $E$. Then 
$$
(mn +l) D\sim  mB +mlE +lD \sim mB + lA,
$$
contradicting 
$E\subset\SBs(D)$.
\qed

\medskip
A curve 
$E\subset X$ is called \emph{negative definite} if the intersection matrix 
$(E_i\cdot E_j)$ is negative definite, where 
$E_i\subset E$ are the irreducible components. Here we use Mumford's rational
intersection numbers \cite{Mumford 1961}.

\begin{proposition}
\label{bijective negative}
Let 
$D\in Z^1(X)$ be a Weil divisor, and 
$E\subset X$   a negative definite curve. If 
$D\cdot E_i=0$ for all irreducible components
$E_i\subset E$, then the canonical map 
$H^0(D)\ra H^0(D+E)$ is bijective.
\end{proposition}

\proof
We may assume that 
$E$ is reduced and have to check that the inclusion
$H^0(D)\ra H^0(D+\sum\mu_iE_i)$ is surjective for all integers
$\mu_i\geq 0$.
Suppose there is a linear equivalence
$D+\sum\mu_iE_i\sim A+\sum\lambda_iE_i$ for some curve 
$A\subset X$ not containing any 
$E_i$, and certain integers 
$\lambda_i\geq 0$. Since 
$E$ is negative definite, there is a unique
$\QQ$-divisor 
$\sum\gamma_iE_i$ with
$$
\sum\gamma_iE_i\cdot  E_j=A\cdot E_j \quad\text{ for all $E_j\subset E$}.
$$
By \cite{Giraud 1982}, Equation 7, we have 
$\gamma_i\leq 0$ because 
$A\cdot E_j\geq 0$. Consequently, 
$\sum\mu_iE_i\cdot E_j = \sum(\gamma_i+\lambda_i)E_j\cdot E_j$, therefore 
$\mu_i\leq\lambda_i$. In other words, 
$\sum\mu_iE\subset A+\sum\lambda_iE_i$, hence 
$H^0(D)\ra H^0(D+E)$ is bijective.
\qed

\begin{remark}
A curve $E\subset X$ is called \emph{contractible} if there is a proper 
birational 
morphism 
$f:X\ra Y$ of proper normal algebraic 
surfaces so that $X-R$ is the isomorphism 
locus.
The preceding results imply that the graded ring $R(D)$ does not change under 
certain contractions. Namely, if either 
$E\subset\SBs(D)$, or $D\cdot E_i=0 $ for 
all irreducible components $E_i\subset E$, then the canonical injection 
$R(D)\subset R(f_*(D))$ is bijective.
\end{remark}

%===========================================================
\section{Rational maps defined by Weil divisors}
\label{raional maps defined by Weil divisors}

Fix a proper normal algebraic surface $X$ and a Weil divisor $D\in Z^1(X)$. 
The 
task now is to describe, in geometric terms, the rational map 
$r_D:X\dashrightarrow P(D)$. The following decomposition is useful for this. 
For each $n\geq 0$ so that $H^0(nD)\neq 0$, write 
$$
nD\sim M_n+ F_n,
$$
where $F_n\subset X$ is the \emph{fixed part}, and $M_n=nD-F_n $ is the 
\emph{movable part} of $nD$. Note that $H^0(M_n)=H^0(nD)$ and that 
$M_n\cdot C\geq 0$  for all curves $C\subset X$. Furthermore, decompose 
$$
F_n = F_n'+ F_n'',
$$
where $F_n'\subset F_n$ is the part of $F_n$ consisting of all
connected components $C\subset F_n$ 
with $C\cdot M_n>0$.
For $n$ sufficiently divisible, the support of $F_n$ is the 1-dimensional part 
of 
$\SBs(D) $.

\begin{proposition}
\label{maximal contractible}
Suppose  that $P(D)$ is a surface and 
assume $\Bs(nD)=\SBs(D)$. Then there is a  
maximal reduced curve 
$R\subset X$ with $M_n\cdot R=0$ and $F_n'\cap R=\emptyset$. Furthermore,  $R$ 
is negative definite and  contractible.
\end{proposition}

\proof
Fix a connected component 
$C\subset F_n'$ and decompose $C=C_1+\ldots +C_r$ into 
irreducible components.
Rearranging indices and allowing repetitions, we may assume that $M_n\cdot 
C_1>0$ and that the intersections $C_i\cap C_{i+1}$ are 0-dimensional.
Choose a rational $\lambda_1>0$ so that $A_1=\lambda_1C_1$ satisfies 
$(M_n+A_1)\cdot C_1>0$.
Inductively, define $A_i=A_{i-1}+\lambda_iC_i$ for some rational number 
$\lambda_i>0$ so that $(M_n+A_i)\cdot C_j>0$ for $1\leq j\leq i$.
We end up with an effective $\QQ$-divisor $A=A_r$ with support $C$.

Repeating this for the other connected components of $F_n'$, we see that there 
is an 
effective $\QQ$-divisor $A$ with support $F_n'$ 
so that $(M_n+A)\cdot C_j>0$ for 
all irreducible components 
$C_j\subset F_n'$. Then 
$$
(M_n\cdot R=0 \quad\text{and}\quad F_n'\cap R=\emptyset)
\quad\Longleftrightarrow\quad
(M_n+A)\cdot R=0
$$
holds for each curve $R\subset X$. By the Hodge Index Theorem, there is a 
maximal reduced 
curve $R\subset X$ satisfying $(M_n+A)\cdot R=0$, and this curve is negative 
definite.

It remains to check that each connected component $R_i\subset R$ is 
contractible.
Choose a curve $A\sim M_n$. Since $M_n$ is movable, each $R_i$ is either 
disjoint from $A\cup F_n'$ or 
entirely contained in $A$. Now the contraction criterion 
\cite{Schroeer 2000}, Proposition 3.3
applies, and we deduce that $R_i$ is contractible. 
\qed

\medskip
The rational map $r_D:X\dashrightarrow P(D)$ has a 
maximal open subset $\dom(r_D)\subset X$ on which it is definable. This open 
subset 
is  called the \emph{domain of definition}. Clearly, 
$X-\SBs(D)\subset \dom(r_D)$. This inclusion, however, might be strict.
To explain this, let $n>0$ be an integer satisfying $\Bs(nD)=\SBs(D)$.
Write
$$
\SBs(D) = \SBs'(D)\cup\SBs''(D)\cup\SBs^0(D),
$$
where $\SBs^0(D)$ is the 0-dimensional part,  $\SBs'(D)$ is the part 
corresponding to $F_n'$, and $\SBs''(D)$ is the part corresponding to $F_n''$.
The next result tells us that this decomposition
depends only on the Weil divisor $D$,
and not on the integer $n>0$:

\begin{theorem}
\label{domain of definition}
The open subset $X-(\SBs(D)'\cup \SBs^0(D))$ is the domain of 
definition for the  rational map $r_D:X\dashrightarrow P(D)$.
Furthermore, the induced  morphism $r_D:\dom(r_D)\ra P(D)$ is proper.
\end{theorem}

\proof
First suppose that $P(D)$ is a curve. Then $M_n^2=0$ and $F_n'=0$ by 
\cite{Schroeer 2000}, Proposition 6.5.
So $M_n$ is a globally generated Cartier divisor, and $P(D)=P(M_n)$.
Consequently, the rational map $r_D$ is everywhere defined.
Furthermore $\SBs^0(D)=\emptyset$,  and the assertion follows.

Now assume that $P(D)$ is a surface. Choose a curve $C\sim M_n$ and set 
$U=X-(C\cup F_n)$.
Then $P(D)$ is covered by affine open subsets of the form $U^\aff$.
Let $A\subset X$ and $R\subset X$ be the curves from the proof of Proposition 
\ref{maximal contractible}, and set 
$$
V= X-(C\cup F_n')=X-(C\cup A).
$$
By \cite{Schroeer 2000}, Proposition 3.2, 
the morphism $V\ra V^\aff$ is proper, 
hence its 
exceptional curve is $R\cap V$. Consequently, $U^\aff=V^\aff$, and $P(D)$  is 
nothing but the contraction of $R\subset X-(F_n'\cup \SBs^0(D))$.
This implies that $X-(F_n'\cup \SBs^0(D))$ is the domain of definition for 
the rational map $r_D:X\dashrightarrow P(D)$.
\qed

\medskip
The preceding proof shows more, namely:

\begin{corollary}
\label{rational map}
Suppose that $P(D)$ is a surface. Then the curve $R\subset X$ from Proposition 
\ref{maximal contractible} is the exceptional set for the proper morphism 
$\dom(r_D)\ra P(D)$.
\end{corollary}

Finally, we mention that the ring $R(D)$ already lives on the  
algebraic surface $P(D)$. In some sense, this reduced the study of
Weil divisors
on proper   surfaces to the study of $\QQ$-Cartier divisors
on algebraic surfaces.

\begin{proposition}
\label{image divisor}
Suppose that $P(D)$ is a surface, and set
$D'= (r_D)_*(D|_{ \dom(r_D)})$. Then the canonical  map $R(D)\ra R(D')$ is 
bijective.
\end{proposition}

\proof
Given an integer $n>0$,
we have to check that the inclusion 
$H^0(nD)\subset H^0(nD+E)$ is bijective for each 
curve $E\subset X $ supported by $\SBs(D)$.
But this follows from Proposition \ref{bijective support}.
\qed

%===========================================================
\section{Compactification}

Let $D$ be a Weil divisor on a proper normal algebraic surface $X$, and assume 
that 
$P(D)$ is 2-dimensional. This algebraic surfaces is not necessarily proper. 
The task 
now is to construct \emph{compactifications} of 
$P(D)$, that is, proper normal surfaces containing $P(D)$ as an open dense 
subset.

We start with a rather simple compactification.
Let $R\subset X$ be the contractible curve from Proposition 
\ref{maximal contractible}  and $f:X\ra Y$ be its contractions.

\begin{proposition}
\label{simple compactification}
There is an open embedding $P(D)\subset Y$, 
and the complement is isomorphic to 
$\SBs'(D)\cup\SBs^0(D)$.
\end{proposition}

\proof
According to Corollary \ref{rational map}, we have a commutative
diagram
$$
\begin{CD}
\dom(r_D) @>>> X\\
@VVV @VVV\\
P(D) @>>> Y,
\end{CD}
$$
where both vertical arrows are the contractions of $R$.
Hence the lower horizontal  arrow exists and 
 gives the desired open embedding.
By Theorem \ref{domain of definition}, its complement
equals $\SBs'(D)\cup\SBs^0(D)$.
\qed

\medskip
The boundary at infinity for $P(D)\subset Y$ contains the 1-dimensional part 
$\SBs'(X)$, so this compactification might not be minimal. The following is 
useful for the construction of smaller compactifications: 

\begin{lemma}
\label{not ample}
Each 1-dimensional connected component 
of $C\subset \SBs(D)$ contains an irreducible component $C_i$ with $D\cdot 
C_i\leq 
0$.
\end{lemma}

\proof
Suppose to the contrary that $D\cdot C_i>0$ for all irreducible components 
$C_i\subset C$.
Let  $f:Y\ra X$ be a resolution of the singularities contained in $C$, and 
$E\subset Y$ the exceptional divisor. Since $E$ is negative definite, 
there is a   
divisor $A\in Z^1(Y)$ supported by $E$ and ample on $E$. 
For $n>0$ sufficiently 
large, the divisor $D'=A+f^*(nD)$ is effective and ample on $f^{-1}(C)$. 
On the other hand, since $nD=f_*(D')$, the preimage $f^{-1}(C)$ contains a 
connected 
component of 
$\SBs(D')$.
This contradicts the Fujita--Zariski Theorem 
(\cite{Fujita 1983}, Theorem 1.19).
\qed

\begin{proposition}
\label{boundary negative}
The 1-dimensional part of $\SBs(D)$ is  negative definite.
\end{proposition}

\proof
We saw in Proposition \ref{maximal contractible} that $\SBs''(D)$  is negative 
definite, so it remains to see
that each connected component $C\subset \SBs'(D)$ is negative definite.
To do so, choose $n>0$ so that $\Bs(nD)=\SBs(D)$, decompose 
$C=C_1+\ldots+C_r$ 
into irreducible components, and let $B\subset F_n'$ be the part disjoint to
$C$.

First, assume that some linear combination $A=\sum\lambda_iC_i$ has $A^2>0$.
As in the proof of \cite{Schroeer 2000}, Proposition 3.2, we may assume that 
$\lambda_i>0$ and 
$A\cdot C_i>0$ for $1\leq i\leq r$.
Replacing $A$ by a suitable multiple, we obtain $F_n'\subset A+B$. By Lemma 
\ref{not ample}, 
there is some component $C_i\subset C$ not contained 
in the stable base locus of 
$M_n+ A+B$. On the other hand, the inclusions
$$
H^0(mD)\supset H^0(m(M_n+F'_n))\subset H^0(m(M_n+A+B))
$$ 
are equalities by Proposition \ref{bijective support}, contradiction.

Second, assume that the intersection matrix $(C_i\cdot C_j)$ is 
\emph{negative semidefinite}. By the Hodge Theorem, its radical  has rank 1.
 This implies that 
there 
is a linear combination
$A=\sum\lambda_iC_i$ with $A\cdot C_i=0$ and $\lambda_i>0$ for all $i$. Up to 
multiples, such a divisor is unique.
Rearranging indices, we may assume $M_n\cdot C_1>0$. The remaining curve 
$C_2+\ldots+C_r$ is negative definite, so there is a linear combination 
$E=\sum\mu_iC_i$ with $E\cdot C_i>0$ for $2\leq i\leq r$. Enlarging $n$ if 
necessary, 
we have $(M_n+E)\cdot C_i>0$ for   $1\leq i\leq r$. 
Furthermore, for some $k>0$, 
the 
divisor $M_n+E+kA+B$ is effective and contains  $F_m'$. 
We conclude the proof as 
in the 
first case.
\qed

\medskip
This implies that the homogeneous spectrum $P(D)=\Proj(R(D))$ 
is not too far from
being proper:

\begin{corollary}
\label{not quasiaffine}
The   algebraic surface $P(D)$ is not quasiaffine.
\end{corollary}

\proof
Suppose to the contrary that $P(D)$ is quasiaffine, and consider the open 
embedding 
$P(D)\subset Y$ from Proposition \ref{simple compactification}.
By Proposition \ref{boundary negative}, the 1-dimensional part 
of the boundary  
is negative definite. On the other hand, \cite{Schroeer 2000}, Proposition 3.2 
ensures that the 
boundary
supports a Weil divisor $C\in Z^1(Y)$ with $C^2>0$, contradiction.
\qed

\medskip
A negative definite curve is not necessarily contractible. However, it can be 
contracted in the category of \emph{algebraic spaces}. 
Roughly speaking, an algebraic space is an object that admits an \'etale 
covering by 
a scheme \cite{Laumon; Moret-Bailly 2000}. Over the complex numbers,    
2-dimensional algebraic spaces 
correspond to \emph{Moishezon surfaces}.

\begin{theorem}
\label{compactification}
There is a proper normal 2-dimensional algebraic space $\overline{P}(D)$ 
containing 
$P(D)$ as an open dense subscheme with 0-dimensional complement.
\end{theorem}

\proof
Since  the 1-dimensional part   $E\subset \SBs(D)$ is negative definite, a 
result of 
Artin (\cite{Artin 1970}, Corollary 6.12) implies that  there is a 
contraction $f:X\ra \overline{P}(D)$ of 
$E\subset X$  in the category of algebraic spaces.
In light of Proposition \ref{rational map}, 
we obtain the desired open embedding 
$P(D)\subset \overline{P}(D)$.
\qed

\begin{remark}
We have constructed a commutative diagram 
$$
\begin{CD}
\dom(r_D) @>>> X\\
@Vr_D VV @VV f V\\
P(D) @>>> \overline{P}(D).
\end{CD}
$$
The boundary at infinity $\overline{P}(D)-P(D)$ comprises the image of 
$\SBs'(D)\cup\SBs^0(D)$. 
Clearly, the points $b\in \overline{P}(D)$ corresponding 
to 
$\SBs^0(D)$ are \emph{scheme-like} \cite{Knutson 1971}, page 131. In contrast, 
the points 
corresponding 
to 
$\SBs'(D)$ might be nonscheme-like. The algebraic space $\overline{P}(D)$ is a 
scheme 
if and only if the negative definite curve $\SBs'(D)$ is contractible.
\end{remark}

Let us collect some properties of the compactification.

\begin{proposition}
\label{Q-Cartier locus}
Suppose that the algebraic space $\overline{P}(D)$ is a scheme. Then the open 
subset 
$P(D)\subset\overline{P}(D) $ is the $\QQ$-Cartier locus for  the Weil
divisor $f_*(D)$.
\end{proposition}

\proof
Clearly, $f_*(D)$ is $\QQ$-Cartier on $P(D)$. 
According to Proposition \ref{not 
Q-Cartier}, it 
is 
not $\QQ$-Cartier on each point of $\SBs^0(D)$.
Suppose there is a connected component $C\subset \SBs'(D)$ so that $f_*(D)$ is 
$\QQ$-Cartier on the image $y=f(C)$.
For $n>0$ sufficiently divisible, $f_*(nD)=f_*(M_n)$ is movable. Again by 
Proposition 
\ref{not Q-Cartier}, we find a curve $E\sim f_*(nD)$ disjoint from $y$. Since 
$R(D)=R(f_*(D))$, 
the corresponding curve $E'\sim nD$ contains $C$ as a connected component, 
contradicting the definition of $\SBs'(D)$.
\qed

\medskip
In light of Nakai's ampleness criterion, we call a Weil divisor 
$D\in Z^1(X)$ \emph{numerically ample} if 
$D^2>0$ and 
$D\cdot C>0$ for all curves 
$C\subset X$.

\begin{proposition}
\label{numerically ample}
Suppose that the algebraic space $\overline{P}(D)$ is a scheme. Then  the Weil
divisor $f_*(D)$ 
is 
numerically ample.
\end{proposition}

\proof
Choose  $n>0$ so that $\Bs(nD)=\SBs(D)$. Then 
$f_*(nD) =f_*(M_n)$ is movable, so 
$f_*(D)\cdot C\geq 0$ for all integral curves $C\subset \overline{P}(D)$. 
Suppose 
$f_*(D)\cdot C= 0$. Then 
$C$ is disjoint from the  non-$\QQ$-Cartier locus of $f_*(D)$. 
Consequently, the effective $\QQ$-divisor $f^*(C)$ is disjoint from $F_n'$.
The projection formula gives 
$$
 M_n\cdot f^*(C) = f_*(nD)\cdot C=0.
$$
By Proposition \ref{maximal contractible}, 
the curve $ f^*(C)$ is contracted by $f:X\ra 
\overline{P}(D)$, 
which is 
absurd. Therefore, $ f_*(D)\cdot C> 0$ for all curves $C$. Then $f_*(D)^2>0$ 
holds as 
well, because
$f_*(nD)$ is effective.
\qed

\begin{remark}
The preceding two Propositions hold true without the assumption
that $\P(D)$ is a scheme. I leave it to the interested reader to
formulate the corresponding results.
\end{remark}

%===========================================================
\section{Multiplicities inside the stable base locus}

We keep our proper normal algebraic surface $X$ 
and  Weil divisor $D\in Z^1(X)$.
Decompose $nD=M_n+F_n$ and $F_n=F_n'+F_n''$ as in Section 
\ref{raional maps defined by Weil divisors}. 
Zariski \cite{Zariski 1962} pointed out that  multiplicities   in $F_n$ play  
role for the structure of  the graded ring 
$R(D)$. First, we start with the part $\SBs''(D)$ of the stable base locus.

\begin{proposition}
\label{unbounded multiplicities}
Suppose $\Bs(nD)=\SBs(D)$. Then $F_{mn}''=mF_n''$ for all $m>0$. 
In particular, 
the multiplicities in $F_{mn}''$ are unbounded.
\end{proposition}

\proof 
Since $F''_n=\SBs(D)$ is negative definite and $mnD\cdot C=mF_n''\cdot C$
 for all curves $C\subset \SBs''(D)$, we infer $F''_{mn}=mF''_n$.
\qed

\medskip
Next, we turn to the part $\SBs'(D)$ of the stable base locus.
Suppose $\Bs(nD)=\SBs(D)$. Decompose $\SBs'(D)=C_1+\ldots+C_r$ into integral 
components.
For each such component, define
$$
\lambda_i=\inf \left\{ \mult_{C_i}(F'_n-A)\mid \text{$A$ as below}\right\},
$$
where $A$ runs through all effective $\QQ$-divisors contained in $F_n'$ 
such that  that the condition
$(M_n+A)\cdot C_j> 0$ holds for $1\leq j\leq r$.

\begin{proposition}
\label{bounded multiplicities}
With the preceding notation,
 $\liminf_{m\ra\infty} \mult_{C_i}(F'_{mn}/m)\leq \lambda_i$.
\end{proposition}

\proof
Let $A$ be an effective $\QQ$-divisor as above. 
By Lemma \ref{not ample}, the curve  
$C_i$ is not contained in $ \SBs(M_n+A)$. 
Now the decomposition into effective summands
$mnD=m(M_n+A) + m(F_n-A)$  gives 
$$
\mult_{C_i}(F_{mn}/m) \leq \mult_{C_i}(F_n-A)
$$
for $m$ sufficiently divisible, hence the assertion.
\qed

%===========================================================
\section{The  canonical model}

Fix a proper normal algebraic surface $X$. In this section, we shall apply
the  results of the preceding Sections to the 
\emph{canonical ring} $R(K_X)=\bigoplus H^0(nK_X)$ 
and the 
corresponding \emph{canonical model} 
$$
P(K_X)=\Proj(R(K_X)).
$$
Note that our canonical 
model is defined on the surface $X$ itself, and not on a resolution of 
singularities.
The algebraic space $\overline{P}(K_X)$ constructed in Proposition 
\ref{compactification} is 
called the \emph{compactified canonical model}. 
The task now is to determine whether or not $\P(D)$ is a scheme.

To do this, the following notions are useful.
Two Weil divisors 
$A,B\in Z^1(X)$ on a proper normal algebraic surfaces $X$ are called 
\emph{numerically equivalent} if 
$A\cdot C=B\cdot C$ for all curves 
$C\subset X$. Given a subset $S\subset X$, we say that $X$ 
is \emph{numerically 
$\QQ$-factorial} with respect to $S$ if each Weil divisor on 
$X$ is numerically equivalent to a $\QQ$-divisor 
that is $\QQ$-Cartier near $S$.
If this holds for $S=X$, we call $X$ numerically $\QQ$-factorial. 
For example, $\QQ$-factorial  surfaces
or surface of geometric genus $p_g=0$ are 
numerically $\QQ$-factorial.

\begin{theorem}
\label{canonical scheme}
If $X$ is numerically $\QQ$-factorial with respect to $\SBs'(K_X)$, then the 
proper algebraic
space $\overline{P}(K_X)$ is a scheme.
\end{theorem}

\proof
We have to check that the negative definite curve $C=\SBs'(K_X)$ is 
contractible.
According to Proposition \ref{bijective support}, the   inclusion 
$H^0(K_K)\subset 
H^0(K_X+C)$ is bijective.
But now the contraction criterion \cite{Schroeer 2000}, Theorem 5.1 applies,
and we deduce that $C\subset X$ is contractible.
\qed

\medskip
To obtain examples, we shall relate numerically 
$\QQ$-factorial surfaces   to surfaces with rational singularities.
 
\begin{proposition}
\label{log terminal discrepancies}
A proper normal algebraic surface is 
$\QQ$-factorial on the locus of rational singularities.
\end{proposition}
 
\proof 
Over an algebraically closed ground field, this easily follows
from \cite{Artin 1962}, Theorem 1.7.
The following argument using group schemes works over
arbitrary ground fields $k$.
Let $f:Y\ra X$  be a resolution of the rational singularities,
$E\subset Y$ the reduced exceptional curve, and 
$\foY\subset Y$ the corresponding formal completion.
Then $H^1(E,\O_{nE})=0$ for all $n>0$. 
Hence the group schemes $\Pic^0_{nE/k}$ vanish, and  $\Pic^0(\foY)=0$.
Given a Weil divisor $D\in Z^1(X)$, there is an integer $n>0$ so that 
the $\QQ$-divisor $C=f^*(nD)$ has integral coefficients. 
Beeing numerically trivial, the
formal line bundle $\O_\foY(C)$ is trivial.
Hence $f_*\O_Y(C)$ is invertible and $f_*(C)=nD$ is Cartier 
on the locus of rational singularities.
\qed
 
\medskip
Together with Theorem \ref{canonical scheme}, this gives 
a condition for schematicity.

\begin{corollary}
\label{log terminal proper}
Suppose $X$ has rational singularities
along $\SBs'(K_X)$. Then the algebraic space
$\overline{P}(K_X)$ is a proper scheme.
\end{corollary}

We also have a condition for projectivity:

\begin{corollary}
\label{log terminal projective}
Suppose $X$ has rational singularities
along $\SBs'(K_X)\cup\SBs^0(K_X)$. Then the algebraic space
$\overline{P}(K_X)$ is a projective scheme.
\end{corollary}

\proof
By Proposition 
\ref{numerically ample}, the canonical class $K$ of $\P(K_X)$ 
is numerically ample.
The locus where $K$  is not $\QQ$-Cartier is nothing but
the image of $\SBs'(K_X)\cup\SBs^0(K_X)$. Hence $K$ is 
numerically equivalent to a $\QQ$-Cartier divisor,
so the compactified canonical model is projective.
\qed

\medskip
The canonical model of a regular surfaces  is not necessarily regular. Rather, 
it has rational Gorenstein singularities. However,   
canonical models of surfaces with   
rational Gorenstein singularities have
rational Gorenstein singularities, so this class of surfaces is closed under 
passing to
minimal models. Here is a similar result.
 
\begin{proposition}
\label{numerically q-factorial}
If  the surface
$X$ is numerically 
$\QQ$-factorial, then  the scheme
$\P(K_X)$ is numerically 
$\QQ$-factorial as well. 
\end{proposition}

\proof
I claim that the negative definite curve 
$C=\SBs'(K_X)\cup\SBs''(K_X)$ is contractible.
Indeed: By Proposition \ref{bijective support}, 
the curve $C$ is contained in the fixed curve of $K_X+mC$ for all $m>0$. 
Furthermore, we assume that
$X$ is numerically $\QQ$-factorial. Hence \cite{Schroeer 2000}, Theorem 5.1 
ensures that $C\subset X$ is contractible.

Next, we check that the corresponding contraction $h:X\ra Y$ 
yields a numerically $\QQ$-factorial surface. To see this, consider for each 
$m>0$ the exact sequence 
$$
H^1(X,\O_{X}) \lra H^1(mF,\O_{mF})  \lra 
H^2(X,\O_{X}(-mF))  \lra H^2(X,\O_{X}).
$$
The map on the right is Serre dual to 
$H^0(K_X)\ra H^0(K_X+mF)$, which is surjective by 
Proposition \ref{bijective support}.
Hence the map on the left is surjective. By Grothendieck's Existence Theorem, 
the map 
$H^1(X,\O_{X})\ra H^1(\foX,\O_{\foX})$ is surjective as well, where 
$\foX=X_{/F}$ is the formal completion.  Now 
\cite{Schroeer 2000}, Proposition 4.2 
ensures that  
$Y$ is numerically 
$\QQ$-factorial.

Recall that the contraction
$f:X\ra\P(K_X)$ admits a factorization $g:Y\ra \P(K_X)$.
Furthermore, $K_Y$ has trivial intersection number on each irreducible
component of the exceptional curve $E\subset Y$.
Now Proposition \ref{bijective negative} ensures that $H^0(K_Y)=H^0(K_Y+mE)$,
 and we 
deduce as 
in the preceding paragraph
that the compactified canonical model $\P(K_X)$ is numerically $\QQ$-factorial.
\qed

\begin{remark}
The numerical criterion for ampleness implies that numerically $\QQ$-factorial 
surfaces are projective. Therefore, their compactified canonical model
is projective as well. More precisely, a multiple of the canonical class
of $\P(K_X)$ deforms to an ample invertible sheaf.
\end{remark}

Let me record the following special case of Propositions
\ref{Q-Cartier locus} and \ref{numerically ample}.

\begin{proposition}
\label{locus}
Suppose that $\P(K_X)$ is a scheme. Then the open subset 
$P(K_X)$ is its  
$\QQ$-Gorenstein locus, and the canonical class of $\P(K_X)$ is numerically 
ample.
\end{proposition}

A Weil divisor $A\in Z^1(X)$ is called \emph{nef} 
if $A\cdot R\geq 0$ holds for 
all curves $R\subset C$.
We have the following characterization for the compactified canonical model to 
be a scheme.

\begin{theorem}
The algebraic space $\P(K_X)$ is a scheme if and only if for each 
connected component $R\subset \SBs'(K_X)$, there is a nef Weil divisor $A\in 
Z^1(X)$ that is Cartier near $R$, 
so that for each integral curve $C\subset X$, the condition $A\cdot 
C=0$ holds if and only $C\subset R$.
\end{theorem}

\proof
You easily check that the condition is necessary.
For sufficiency, we shall apply the characterization of contractible curves
in \cite{Schroeer 2000}, Theorem 3.4 to  each connected component 
$R\subset\SBs'(D)$.
Let $R$ be as in the Theorem. Fix an integer $m>0$ and consider the exact 
sequence
$$
H^1(X,\O_X)\lra H^1(mR,\O_{mR})\lra H^2(X,\O_X(-mR))\lra H^2(X,\O_X).
$$
The map on the right is Serre dual to $H^0(K_X)\ra H^0(K_X+mR)$. The latter is 
surjective according to Proposition \ref{bijective support}. 
Consequently, $H^1(X,\O_X)\ra 
H^1(mR,\O_{mR})$ is surjective, so 
$\Pic^0(X)\ra \Pic^0(mR)$ is surjective up to 
torsion. So, for some $n>0$, the invertible sheaf $\O_{mR}(nA)$ is the 
restriction of some numerically trivial invertible $\O_X$-module. Thus the 
conditions of the before mentioned characterization of contractible curves 
applies, and we conclude that  $C$ is contractible.
\qed

\begin{remark}
One might say that a proper normal algebraic surface 
$X$ with rational singularities has \emph{two}  canonical models: 
First the   
canonical model 
$P(K_{Y})=\P(K_Y)$ in the sense of Mori theory defined using a resolution of 
singularities 
$Y\ra X$, and second the compactification 
$P(K_{X})\subset\P(K_X)$ defined on the normal surface $X$ itself.
I wish to know what happens in higher dimensions.
\end{remark}

\begin{question}
Does there exist a proper normal surface whose compactified canonical model
is not a scheme?
\end{question}

%===========================================================

\end{document}